\documentclass[twoside,11pt]{amsart}
\usepackage[latin1]{inputenc}
\usepackage{amssymb,amsfonts,amsmath,amsthm,amscd, epsfig}

\usepackage[all]{xy}

\title{On the Galois Map for Groupoid Actions}
\author[Paques and Tamusiunas]{Antonio Paques and Tha\'{i}sa Tamusiunas}
\address{Instituto de Matem\'{a}tica, Universidade Federal do Rio Grande do Sul,  Av. Bento Gon\c{c}alves, 9500, 91509-900. Porto Alegre-RS, Brazil}
\email{paques@mat.ufrgs.br, trtamusiunas@yahoo.com.br}
\date{}

\newtheorem{teo}{Theorem}[section]

\newtheorem{lema}[teo]{Lemma}

\newtheorem{prop}[teo]{Proposition}

\newtheorem{obs}[teo]{Remark}

\newcommand{\vu}{\vspace{.1cm}}
\newcommand{\vd}{\vspace{.2cm}}

\newcommand{\G}{\mathcal{G}}

\begin{document}

\maketitle

\begin{abstract}
Some conditions for the Galois map to be injective are given in the groupoid acting on a noncommutative ring context. In the particular case in which the Galois extension is a central Galois algebra, it is given a complete characterization of that kind of extension with Galois map bijective.
\end{abstract}

\vspace{0,5 cm}


\section{Introduction}

Since the work published in 1926 by H. Brandt in \cite{brandt}, the study of the theory of groupoids  is growing and gaining developments in different areas of mathematics and even physics, as shown in studies involving algebraic topology, noncommutative geometry, Lie groupoids, partial actions and Galois theory (\cite{bagio}, \cite{brown}, \cite{brown2}, \cite{brown3}, \cite{corta}, \cite{moe}, \cite{pinedo}, \cite{pata}).

\vu

Actions of groupoids on rings were introduced in the literature by D. Bagio and the first author in \cite{bagio}. Following them, an action of a groupoid $\G$ on a ring $R$ is a pair $\beta = (\{E_g\}_{g \in \G}, \{\beta_g\}_{g \in \G})$, where each $E_g$ is an ideal (not necessarily nonzero) of $R$ and each $\beta_g$ is an isomorphism of ideals from $E_{g^{-1}}$ onto $E_g$ satisfying some appropriate conditions of compatibility. In order to attend our purposes we will also assume along all this manuscript that such ideals $E_g$ are all unital, with the corresponding identity element denoted by $1_g$. We also recall from \cite{bagio} that any action $\beta$ of $\G$ on $R$ induces by restriction the action $\beta_H=(\{E_h\}_{h \in H}, \{\beta_h\}_{h \in H})$ of $H$ on $R$, for any subgroupoid $H$ of $\G$. Notice, in particular, that $\beta_{\G}=\beta$.

\vu

Any action of a groupoid $\G$ on a ring $R$ gives rise to a natural map from the set of the subgroupoids of $\G$ into the set of  subrings of $R$ via a  standard and well known procedure that in the context of groupoid actions on rings has the following interpretation: to each subgroupoid $H$ of $\G$ is associated the subring, denoted by $R^{\beta_H}$, of all the elements $r$ of $R$ satisfying the condition $\beta_h(r1_{h^{-1}})=r1_h$\, for all $h\in H$, the so called \emph{invariants of $R$  by the action $\beta_H$ of $H$ on $R$.} Such a map is clearly an inverting-inclusion map and by  consequence the subring $R^\beta$ of all the invariants of $R$ by $\beta$ is a subring of $R^{\beta_H}$, for any subgroupoid $H$ of $\G$. In the case that $R$ is a Galois extension of $R^\beta$ this map is called the \emph{Galois map for the action $\beta$ of $\G$ on $R$} and denoted by $\theta_\beta$ or simply $\theta$ if no confusion arises.

 \vu

It was already proved by the authors in \cite{paqtamII} that in the commutative context there exists a bijective correspondence between the wide subgroupoids of $\G$ and the $R^\beta$-subalgebras of $R$ that are separable and $\beta$-strong. This earlier result is an extension to the setting of groupoid actions of the classical and celebrated Galois correspondence theorem due to S.U. Chase, D.K. Harrison and A. Rosenberg in \cite{chase} for group actions.

\vu

Furthermore, still for group actions but in the noncommutative context, G. Szeto and L. Xue proved in \cite{xue} the bijectivity of the Galois map in the context of  central Galois algebras and later in \cite{szeto} they presented necessary and sufficient conditions for the injectivity of the Galois map in the more general context of Galois extensions. The main tool that plays an important role in these works is the specifical set of  all $C(R)$-modules $J_g = \{r \in R \mid rg(x) = xr$, for all $x \in R\}$, where $C(R)$ denotes  the center of $R$ and each $g$ an element of the group acting on $R$.

 \vu

For a groupoid $\G$ acting on a ring $R$ via an action $\beta = (\{E_g\}_{g \in \G}, \{\beta_g\}_{g \in \G})$,  we define $J_g$ as being the set of the elements $r$ of $E_g$ such that $r\beta_g(x1_{g^{-1}}) = xr$, for all $x \in R$.

\vu

This paper is organized as follows. In the next section, essentially, we recall from the literature a few of definitions necessary for all what follows, concerning to groupoid actions, Galois extensions, Galois algebras and separability, mainly Hirata separability. Section 3 is dedicated to discuss conditions necessary and sufficient, in terms of $J_g$, $g \in \G$, under which the Galois map for groupoid actions on rings can be injective (see Theorems
\ref{thetainj}, \ref{hinjetiva}, \ref{problema} and \ref{gammatheta}). In Section 4 we present a characterization of the central $\beta$-Galois algebras for which the Galois map is bijective. Actually, we prove that a central $\beta$-Galois extension $R$ of $R^\beta$ satisfies the Galois correspondence theorem if and only if for each separable $R^{\beta}$-subalgebra $S$ of $R$, the commutator $V_R(S)=\{r\in R\ |\ rs=sr,\, \text{for all}\, s\in S\}$ of $S$ in $R$  is equal to  the direct sum of all $J_g$, with $g\in H_S$, where $H_S = \{g \in \G\ |\ \beta_g(s1_{g^{-1}}) = s1_g,\, \text{for all}\, s \in S\}$ (see Theorem \ref{teolast}). All these last results generalize the work developed by Szeto and Xue in \cite{xue} and \cite{szeto}.

\vu

Throughout, unless otherwise specified, rings (hence, also algebras) are associative and unital.


\section{Prerequisites and Basic Results}

Groupoids are usually presented as small categories whose morphisms are all invertible. However, they can also be seen as a generalization of groups. We will adopt here such a version which appears, for instance, in \cite{lawson}. A \emph{groupoid} is a nonempty set $\G$, equipped with a partially defined binary operation, which we will denote by concatenation, that satisfies the associative law (whenever this makes sense) and the condition that every element $g \in \G$ has an inverse $g^{-1}$ and a left and a right  identity, respectively denoted by $r(g)$ and $d(g)$.

\vu

For all $g, h \in \G$, we write $\exists gh$ whenever the product $gh$ is defined. As it was observed, for instance in \cite{pata}, it immediately follows from the above definition that for all $g, h \in \G$, $\exists gh$ if and only if $d(g) = r(h)$ and, in this case, $d(gh) = d(h)$ and $r(gh) = r(g)$. We will be denoted by $\G^2$ the subset of the pairs $(g, h) \in \G \times \G$ such that $d(g) = r(h)$.

An element $e \in \G$ is called an \emph{identity} of $\G$ if $e = d(g) = r(g^{-1})$, for some $g \in \G$. It will be denoted by $\G_0$ the set of all identities of $\G$ and by $G_e$ the set of all $g \in \G$ such that $d(g) = r(g) = e$. Clearly, $\G_e$ is a group, called the \emph{principal (or isotropic) group associated to e}.

\vu

Given  a groupoid $\G$ and  a nonempty subset $H$ of $\G$, we say that $H$ is a \emph{subgroupoid} of $\G$ if the following conditions are satisfied:
\begin{itemize}
\item [(i)] For all $g, h \in H$, if $\exists gh$ then $gh \in H$;
    \item [(ii)] For all $g \in \G$, if $g \in H$ then $g^{-1} \in H$.
\end{itemize}
If $H_0 = \G_0$, we say that $H$ is \emph{wide} (\cite{pata}, p. 85).

\vu

Consider $R$ an algebra over a commutative ring $K$. Recalling from \cite{bagio}, an \emph{action} \emph{of} $\G$ \emph{over} $R$ is a pair $$\beta = (\{E_g\}_{g \in \G}, \{\beta_g\}_{g \in \G})$$ where for each $g \in \G$, $E_g = E_{r(g)}$ is an ideal of $R$ and $\beta_g: E_{g^{-1}} \to E_g$ is an isomorphism of $K$-algebras satisfying the following conditions: \begin{itemize}
\item [(i)] $\beta_e$ is the identity map $Id_{E_e}$ of $E_e$ for all $e \in \G_0$;
    \item [(ii)] $\beta_g(\beta_h(r)) = \beta_{gh}(r)$ for all $(g, h) \in \G^2$ and for all $r \in E_{h^{-1}} = E_{(gh)^{-1}}$.\end{itemize}
If each $E_g$, $g\in\G$,  is a unital algebra (with identity element denoted by $1_g$) we say that the action $\beta$ is \emph{unital}.

From now on, along all the text, any action of a groupoid $\G$ on an algebra $R$ will be assumed unital.

\begin{lema}\cite[Lemma 2.1]{paqtamII}\label{ahat} Consider $\beta = (\{E_g\}_{g \in \G}, \{\beta_g\}_{g \in \G})$ a action of a groupoid $\G$ on an algebra $R$.  For any subalgebra $T$ of $R$, define $$H_T = \{g \in \G \mid \beta_g(t1_{g^{-1}}) = t1_g,\, \text{for all}\,\  t\in T\}.$$ Then $H_T$ is a subgroupoid of $\G$.
\end{lema}

\vu

Following \cite{flores}, the \emph{skew groupoid ring} $R \star_{\beta} \G$, corresponding to an action $\beta$ of a groupoid $\G$ on an algebra $R$, is defined as the direct sum $$R \star_{\beta} \G = \bigoplus_{g \in \G}E_gu_g,$$ (where each $u_g$ is a placeholder for the g-th component) with the usual addition and the multiplication induced by the rule
$$
(xu_g)(yu_h) = \left\{
\begin{array}{lll}
x\beta_g(y1_{g^{-1}})u_{gh},& \mbox{if} \quad (g, h) \in G^2\\
0, & \mbox{otherwise,}
\end{array}
\right.
$$ for all $g, h \in \G, x \in E_g$ and $y \in E_h$. It is straightforward to check that $R \star_{\beta} \G$ is associative and, if $\G_0$ is finite, also unital, with identity element given by $1_{R \star_{\beta} \G}=\sum_{e\in\G_0}1_e\delta_e$.

\vd

For any action $\beta$ of  a groupoid $\G$ on an algebra $R$ we will denote by $$R^{\beta}=\{r\in R\ |\ \beta_g(r1_{g^{-1}})=r1_g,\,\ \text{for all}\,\ g\in\G\}$$ the subalgebra of $R$ of the elements which are invariant under $\beta$.

\vu

We say that $R$ is a $\beta$-\emph{Galois extension of} $R^{\beta}$ if $\G$ is finite and there exist elements $x_i, y_i \in R$, $1 \leq i \leq m$, such that $\sum_{1 \leq i \leq m}x_i\beta_g(y_i1_{g^{-1}}) = \delta_{e, g}1_e$ for all $e \in \G_0$ and $g \in \G$ \cite{bagio}. The set $\{x_i, y_i\}_{1\leq i \leq m}$ is called a \emph{Galois coordinate system} of $R$ over $R^{\beta}$.

\vu

Let $R\supseteq S$ be a ring extension. We say that $R$ is \emph{separable over $S$} if there exists an element $z = \sum_{i = 1}^nx_i \otimes_S y_i \in R \otimes_S R$ such that $\sum_{i = 1}^nx_iy_i = 1_R$ and $rz = zr$, for every $r\in R$ \cite{HS}. Although in general such an element $z$ is not an idempotent, it is usually referred  in the literature as an \emph{idempotent of separability} of $R$ over $S$.

\begin{prop}\label{propoprimeira} Suppose that $\G$ is a finite groupoid, $K$ a commutative ring, $R$ a $K$-algebra and $\beta = (\{E_g\}_{g \in \G}, \{\beta_g\}_{g \in \G})$ a unital action of $\G$ on $R$.  Assume also that $R$ is a $\beta$-Galois extension of $R^{\beta}$ and $R = \bigoplus_{e \in \G_0}E_e$. Let $H$ be a subgroupoid of $\G$ and $R_H = \bigoplus_{e \in H_0}E_e$. Then $\beta_H = \{\beta_h: E_{h^{-1}} \to E_h \mid h \in H\}$ is an action of $H$ on $R_H$ and $R_H$ is a $\beta_H$-Galois extension of $(R_H)^{\beta_H}$.

\vu

{\bf Proof:} \em  It follows, with the necessary adaptations to the noncommutative case, by the same arguments used in the proof of \cite[Theorem 4.1]{paqtamII}.
\hfill{$\blacksquare$}
\end{prop}

We say that $R$ is a \emph{$\beta$-Galois algebra of $R^{\beta}$} if $R$ is a $\beta$-Galois extension of $R^{\beta}$ such that $R^{\beta}$ is contained in the center $C(R)$ of $R$. Moreover, $R$ is said to be a \emph{central $\beta$-Galois algebra} if $R$ is a $\beta$-Galois extension of its center $C(R)$ (which means that $C(R) = R^{\beta}$). A ring $R$ is called a \emph{Hirata separable extension of $S$} if $R \otimes_S R$ is isomorphic  as an $R$-bimodule to a direct summand of a finite direct sum of copies of $S$ \cite{sugano}. $R$ is said to be a \emph{Hirata separable $\beta$-Galois extension of $R^{\beta}$} if it is a Hirata separable and $\beta$-Galois extension of $R^{\beta}$. We call $R$ an \emph{Azumaya algebra} if it is a separable extension of its center(\cite{demeyer}).

\begin{obs} \label{hseparavel} {\rm Note that, by what was above defined, every central $\beta$-Galois algebra is an Azumaya algebra. Furthermore, by \cite{sugano}, an Azumaya algebra is a Hirata separable extension. Thus, every central $\beta$-Galois algebra is a Hirata separable extension.}
\end{obs}


\section{$\beta$-Galois Extensions with Injective Galois Map}

From now on, we will adopt the following fixed notations: $K$ a commutative ring, $R$ a $K$-algebra, $\G$ a finite groupoid, $\beta = (\{E_g\}_{g \in \G}, \{\beta_g\}_{g \in \G})$ a unital action of $\G$ on $R$  such that $R = \bigoplus_{e \in \G_0}E_e$, $C(R)$ the center of $R$, $V_{S_2}(S_1) = \{r \in S_2 \mid rs = sr$, for all $s \in S_1\}$ the commutator of $S_1$ in $S_2$, for any subrings $S_1, S_2$ of $R$, $J_g$, $g\in\G$, the $C(R)$-submodule of $R$ defined by $J_g = \{r \in E_g \mid r\beta_g(x1_{g^{-1}}) = xr,\, \text{for all}\,\ x \in R\}$ and $\mathcal{S}_H$  the set of all elements $h$ of $H$ such that $J_h\neq 0$, for any subgroupoid $H$ of $\G$.

We also assume that $R$ is a $\beta$-Galois extension of $R^\beta$ and fix the notation $\theta: H \mapsto R^{\beta_H}$ for the Galois map from the set of the wide subgroupoids of $\G$ into the set of the subalgebras of $R$ including $R^{\beta}$. Notice that if $H$ is a wide subgroupoid of $\G$ then $R_H=R$, hence $R$ is a $\beta_H$-Galois extension of $R^{\beta_H}$, by Proposition \ref{propoprimeira}.

\vu

 We start with the following two lemmas, the first one generalizes \cite[Proposition 1]{kanzaki} and it is crucial for all what follows in this section.

\begin{lema}\label{comutador} $V_R(R^{\beta}) = \bigoplus_{g \in \G}J_g.$

\vu

{\bf Proof:} \em Since $R$ is a $\beta$-Galois extension of $R^{\beta}$, by \cite[Theorem 5.3]{bagio}, the map $j: R \star_{\beta}\G \to End(R)_{R^{\beta}}$ given by $j(\sum_{g \in \G}a_g\delta_g)(x) = \sum_{g \in \G}a_g\beta_g(x1_{g^{-1}})$ is a ring isomorphism.
 Furthermore, notice that each $\lambda\in R$ determines uniquely an element $\lambda_l\in End(R)_{R^{\beta}}$ defined by $\lambda_l(r)=\lambda r$, for all $r\in R$. The set $R_l$ of all $\lambda_l$, $\lambda\in R$, is a subring of $End(R)_{R^{\beta}}$ naturally isomorphic to $R$. Analogously, we also have the right isomorphic copy $R_r$ of $R$ as subring of $End_R(R)$, whose elements $\lambda_r$, $\lambda\in R$, are given by $\lambda_r(x)=x\lambda$, for all $x\in R$.

\vu

In order to get our goal  we will proceed by steps.

\vu

\textbf{Step 1:}  $V_{R \star_{\beta}\G }(R) \simeq V_{End(R)_{R^{\beta}}}(R_l)$ via $j$.

\vu

Firstly note that the map $\psi: R \to R \star_{\beta} G$ given by $\psi(r) = \sum_{e \in G_0}r1_e\delta_e$ is injective, thus we may identify $R$ with its isomorphic image in $R \star_{\beta} G$ via such a map and we have that $$j(r)(x) = j(\sum_{e \in G_0}r1_e\delta_e)(x) = \sum_{e \in G_0}r\beta_e(x1_e) = rx\sum_{e \in G_0}1_e = rx1_R=rx = r_l(x),$$ for all $r, x \in R$. Therefore $j(R) = R_l$ and identifying $R \star_{\beta} G$ with its isomorphic image $End(R)_{R^{\beta}}$ via $j$, the result follows.

\vd

\textbf{Step 2:} $V_{End(R)_{R^{\beta}}}(R_l) = End_{R}(R)_{R^{\beta}}$

\vu

It is immediate  to see that $f \in V_{End(R)_{R^{\beta}}}(R_l)$ if and only if $$f(rx) = f \circ r_l(x) = r_l \circ f(x) = rf(x),$$ for all $r, x \in R$, if and only if $f \in End_{R}(R)_{R^{\beta}}$.

\vd

\textbf{Step 3:} $End_{R}(R)_{R^{\beta}} = (V_R(R^{\beta}))_r.$

\vu

$(\subseteq)$ Consider $f \in End_{R}(R)_{R^{\beta}}$. Then $$f(x) = f(x1_R) = xf(1_R)=f(1_R)_r(x),$$ for all $x\in R$. Hence, we have  in particular that   $$f(1_R)_l(x)=f(1_R)x = f(1_Rx)=f(x)=f(x1_R) = xf(1_R)=f(1_R)_r(x),$$ for all $x \in R^{\beta}$, which implies that $f=f(1_R)_r \in V_R(R^{\beta})_r$.

\vu

$(\supseteq)$ Take $f \in (V_R(R^{\beta}))_r$. Then there exists $a \in V_R(R^{\beta})$ such that $f(y) = ya$ for all $y \in R$. Moreover, $a \in V_R(R^{\beta})$ implies that $ax = xa$ for all $x \in R^{\beta}$. Hence, $$f(ryx) = ryxa = ryax = rf(y)x,$$ for all $r\in R$ and $x\in R^\beta$, which ensures that $f \in End_{R}(R)_{R^{\beta}}$.

\vu

On the other hand,

\vd

\textbf{Step 4:} $V_{R \star_{\beta}\G }(R) = \bigoplus_{g \in \G}J_g\delta_g$.

\vu

It is enough to see that
\[
\begin{array}{cll}
\sum_{g \in G} r_g\delta_g \in V_{R \star_{\beta}\G }(R)& \Longleftrightarrow (\sum_{g \in G} r_g\delta_g)x = x(\sum_{g \in G} r_g\delta_g),\,
\forall x \in R\\
&\Longleftrightarrow \sum_{g \in G}r_g\beta_g(x1_{g^{-1}})\delta_g = x\sum_{g \in G} r_g\delta_g,\,\ \forall x \in R\\ & \Longleftrightarrow  r_g\beta_g(x1_{g^{-1}}) = xr_g,\,\forall g\in \G\,\ \text{and}\,\ x\in R\\
&\Longleftrightarrow r_g \in J_g,\, \forall g \in \G\\
&\Longleftrightarrow \sum_{g \in G} r_g\delta_g \in \bigoplus_{g \in \G}J_g\delta_g.\
\end{array}
\]

\vd

\textbf{Step 5:} $\bigoplus_{g \in G}J_g\delta_g \simeq$ $(\bigoplus_{g \in G} J_g)_r$ via $j$.

\vu

 Take $\sum_{g \in G}r_g\delta_g\in\bigoplus_{g \in G}J_g\delta_g$. Since each $r_g\in J_g$, $g\in\G$, it follows that
\[
\begin{array}{ccl}
j(\sum_{g \in G}r_g\delta_g)(x)& = & \sum_{g \in G}r_g\beta_g(x1_{g^{-1}}),\\
 &= &\sum_{g \in G}xr_g = x\sum_{g \in G}r_g\\
 & = &(\sum_{g \in G}r_g)_r(x),\
\end{array}
 \] for all $x\in R$, hence $j(\sum_{g \in G}r_g\delta_g) = (\sum_{g \in G}r_g)_r$.

\vu

Now it is immediate to see from the above steps that $V_{R \star_{\beta}\G }(R) \simeq (\bigoplus_{g \in G}J_g)_r$ via $j$. Therefore, $V_R(R^{\beta}) = \bigoplus_{g \in G}J_g$, which achieve the proof. \hfill{$\blacksquare$}
\end{lema}

\begin{lema}\label{theta} Let $H$ and $L$ be subgroupoids of $\G$ such that $J_h \neq \{0\}$ and $J_l \neq \{0\}$ for every $h\in H$ and $l\in L$. Consider the Galois map $\theta: H \mapsto R^{\beta_H}$. If $\theta(H) = \theta(L)$, then $H = L$.

\vu

{\bf Proof:} \em Since $\theta(H) = \theta(L)$, we have that $R^{\beta_H} = R^{\beta_L}$. In the sequel we will proceed by steps again. In order to simplify notation, we set $H\vee L=\langle\! H\cup L\!\rangle$ to denote the subgroupoid of $\G$ generated by the elements of $H\cup L$, that is, the intersection of all subgroupoids of $\G$ containing $H\cup L$.

\vu

 \textbf{Step 1: } $R$ is a $\beta_{H\vee L}$-Galois extensions of $R^{\beta_{H\vee L}}$

 \vu

 First of all, it is straightforward to check that $\beta_{H\vee L}=\beta_H\cup\beta_L$, hence the Galois coordinate system of $R$ corresponding to $\beta_{H\vee L}$ is precisely the union of those corresponding to $\beta_H$ and $\beta_L$ respectively, and the claim follows.

 \vu

 \textbf{Step 2: } $R^{\beta_{H\vee L}}= R^{\beta_H} = R^{\beta_L}$.

\vu

 The inclusion $R^{\beta_{H\vee L}} \subseteq R^{\beta_H}$ is obvious, for $H \subseteq H\vee L$. The reverse inclusion is also immediate since $\beta_{hl}(r1_{l^{-1}h^{-1}}) = \beta_h(\beta_l(r1_{l^{-1}})1_{h^{-1}}) = \beta_h(r1_l1_{h^{-1}}) = r1_{lh},$  for all $r \in R^{\beta_H} = R^{\beta_L}$, $h\in H$ and $l\in L$ such that and $d(h) = r(l)$.

 \vu

 \textbf{Step 3: } $\bigoplus_{g\in H\vee L} J_g =$ $\bigoplus_{h \in H} J_h$ $= \bigoplus_{l \in L} J_l$.

 \vu

 It follows from Lemma \ref{comutador}, since $V_R(R^{\beta_{H\vee L}}) = V_R(R^{\beta_H}) = V_R(R^{\beta_L})$.

 \vu

 \textbf{\bf Step 4:} $H=L$.

It follows from Step 3, since $J_g\neq 0$ for all $g\in H\cup L$.  \hfill{$\blacksquare$}
\end{lema}

\vu

As an immediate consequence of Lemma \ref{theta} (hence, also of Lemma \ref{comutador}) we have now the first of the main theorems of this section.

\begin{teo} \label{thetainj}If $J_g \neq 0$ for each $g \in G$, then the Galois map $\theta: H \mapsto R^{\beta_H}$ is injective.

\vu

{\bf Proof:} \em By assumption, $J_g \neq 0$ for each $g \in G$. Thus, by the Lemma \ref{theta}, $\theta$ is injective. \hfill{$\blacksquare$}
\end{teo}

\begin{lema}\label{progerador} If $R$ is a Hirata separable extension of $R^{\beta}$, then $J_hJ_g = J_{gh}$ for any $g,h\in\G$ such that $d(g)=r(h)$. In particular, $J_{g^{-1}}J_g = V_{E_{g}}(R)$ and $J_g \neq 0$ for all $g \in G$.

 \vu

{\bf Proof:} \em By \cite[Definiton 1]{sugano}, to say that $R$ is a Hirata separable extension of $R^{\beta}$ is equivalent to say that, for each $g \in G$, the map $$\upsilon_g: V_R(R^{\beta}) \otimes_{C(R)} M^{R} \to M^{R^{\beta}},\quad \upsilon_g(d \otimes m) = dm,\,\,\text{for all}\, d\in V_R(R^{\beta})\,\text{and}\,\ m \in M,$$  where $M$ is a $(R,R)$-bimodule and $M^{R} = \{m \in M\ |\ xm = mx,\,\, \text{for all}\,\ x \in R\}$, is an isomorphism of $C(R)$-modules.

\vu

Note also that from $R$ we obtain an $(R,R)$-bimodule $R_g$ such that $R_g = R$ as a left $R$-module, and as a right $R$-module the action is given by
$$x\cdot y = x\beta_g(y1_{g^{-1}}),\, \text{for}\, x, y \in R.$$

Hence, $J_g = (R_g)^R$, $V_R(R^{\beta}) = (R_g)^{R^{\beta}}$ and the map $\upsilon_g: V_R(R^{\beta}) \otimes_{C(R)} J_g \to V_R(R^{\beta})$ given by $\upsilon_g(d \otimes d_g) = dd_g$, for all $d \in V_R(R^{\beta})$ and $d_g \in J_g$, is an isomorphism of $C(R)$-modules, which ensures, in particular, that $J_g\neq 0$ for each $g\in\G$.

\vu

Furthermore, for every $g, h \in \G$ such that $d(g) = r(h)$, $x\in R$, $r_h \in J_h$ and $r_g \in J_g$ we have: $$xr_hr_g = r_h\beta_h(x1_{h^{-1}})r_g = r_hr_g\beta_g(\beta_h(x1_{h^{-1}})1_{g^{-1}}) = r_hr_g\beta_{gh}(x1_{h^{-1}g^{-1}}),$$ which means that $J_hJ_g\subseteq J_{gh}$.

\vd

On the other hand, since by Lemma \ref{comutador} $V_R(R^{\beta}) = \bigoplus_{g \in G}J_g$, then the isomorphism of $C(R)$-modules $$V_R(R^{\beta})J_g \simeq V_R(R^{\beta}) \otimes_{C(R)} J_g \simeq V_R(R^{\beta}),$$ above described, implies that $$V_R(R^{\beta}) = \bigoplus_{h \in G}J_hJ_g \subseteq \bigoplus_{r(h)=d(g)}J_{gh}\subseteq \bigoplus_{r(l) = r(g)}J_l\subseteq \bigoplus_{l \in G}J_l = V_R(R^{\beta})$$

\vu

Now, it is easily seen that necessarily $J_hJ_g = J_{gh}$ for all $g,h\in\G$ such that $d(g)=r(h)$ and, in particular, $J_{g^{-1}}J_g = J_{r(g)} =  V_{E_g}(R)$, for all $g\in\G$. The proof is complete. \hfill{$\blacksquare$}
\end{lema}

\vu

As an  application of Lemma \ref{progerador} we have the second main theorem of this section.

\begin{teo}\label{hinjetiva} If $R$ is a Hirata separable $\beta$-Galois extension, or a central $\beta$-Galois algebra, then $\theta$ is injective.

\vu

{\bf Proof:} \em For any Hirata separable $\beta$-Galois extension, we have by Lemma \ref{progerador} that $J_g \neq \{0\}$ for each $g \in G$. So, $\theta$ is injective by Theorem \ref{thetainj}. Since every central $\beta$-Galois algebra is a Hirata separable $\beta$-Galois extension (see Remark \ref{hseparavel}), it follows that $\theta$ is injective also in this case. \hfill{$\blacksquare$}
\end{teo}

For the sequel we introduce the following notations:
\begin{itemize}
\item $\sigma(h) = J_h$ for all $h\in H$

\item $\gamma(H) = \bigoplus_{h \in H}\sigma(h)$,
\end{itemize}
for any subgroupoid $H$ of $\G$. Actually, $\sigma$ denotes a map from $H$ to the set of all $J_g$, $g\in \G$, and $\gamma$ a map from the set of the subgroupoids of $\G$ into the set of all $C(R)$-modules. In the sequel we will see some properties of $\sigma$ and $\gamma$ that provide sufficient conditions for $\theta$ to be injective. We start with the following.

\begin{lema} For any subgroupoid $H$ of $\G$, the restriction of the map $\sigma$ to the set $\mathcal{S}_H$  is injective.

\vu

{\bf Proof:} \em Let $J_g = J_h$ for $g, h \in \mathcal{S}_H$. Then, $E_g \cap E_h \neq 0$ and consequently $r(g) = r(h)$, for $R$ is the direct sum of all
$E_e$, $e\in\G_0$.  Since $J_g = J_h \neq 0 $, there exists $0 \neq r \in E_g = E_h$ such that $r \in J_g = J_h$ and $xr = r\beta_g(x1_{g^{-1}}) = r\beta_h(x1_{h^{-1}})$ for all $x \in R$. Thereby, $r(\beta_g(x1_{g^{-1}}) - \beta_h(x1_{h^{-1}})) = 0$  or, equivalently,
$r\beta_g(x1_{d(g)} - \beta_{g^{-1}h}(x1_{h^{-1}g})) = 0$, for all $x \in R$.  Notice that $\exists g^{-1}h$, since $r(g) = r(h)$. Also observe that as $r \in J_g$, $r\beta_g(x1_{g^{-1}}) = xr$ for all $x \in R$, so $$r\beta_g(x1_{d(g)} - \beta_{g^{-1}h}(x1_{h^{-1}g})) = (x - \beta_{g^{-1}h}(x1_{h^{-1}g}))r = 0. \quad (\ast) $$

Since $R$ is a $\beta$-Galois extension of $R^{\beta}$, there exist $x_i, y_i \in R$, $1 \leq i \leq m$, such that \begin{center}$\sum_{1 \leq i \leq m}x_i\beta_g(y_i1_{g^{-1}}) = \delta_{e, g}1_e$, for all $e \in \G_0$ and $g \in \G$. \end{center} Taking $x$ as $y_i$ in the equation $(\ast)$, we get $(y_i1_{d(g)} - \beta_{g^{-1}h}(y_i1_{h^{-1}g}))r = 0$. Thus, $x_i(y_i - \beta_{g^{-1}h}(y_i1_{h^{-1}g}))r = 0,$ for $1 \leq i \leq m$,  which implies that \begin{center}$(\sum_{i = 1}^m x_iy_i)r = \sum_{i = 1}^m x_i\beta_{g^{-1}h}(y_i1_{h^{-1}g})r.$\end{center} Then,
$$
\begin{array}{ccl}
r & \! = & \! 1_{r(g)}r = (\sum\limits_{1\leq i\leq m} x_i\beta_{r(g)}(y_i1_{r(g)}))r = (\sum\limits_{1\leq i\leq m} x_iy_i1_{r(g)})r  \\
\vu
& \! = & \! (\sum\limits_{1\leq i\leq m} x_iy_i)r = \sum\limits_{1\leq i\leq m} x_i\beta_{g^{-1}h}(y_i1_{h^{-1}g})r\\ &\!=&\! \delta_{e, g^{-1}h}1_er,
\end{array}
$$
for all $e \in G_0$. Since $r \neq 0$, necessarily $g^{-1}h\in\G_0$, and so, in particular, $g^{-1}h = r(g^{-1}h) = r(g^{-1}) = d(g)$. Hence,
$h = r(h)h = r(g)h = gg^{-1}h = gd(g) = g$. The proof is complete.  \hfill{$\blacksquare$}
\end{lema}

By keeping the notations, consider $\gamma: H \mapsto \bigoplus_{h \in H}J_h$.

\begin{lema}\label{igualdadegamma} Let $H$ be a wide subgroupoid of $\G$. Then $\gamma(H)$ is a subalgebra of $R$ over $C(R)$ and $\gamma(H) = V_R(R^{\beta_{\mathcal{S}_H}})$, where $R^{\beta_{\mathcal{S}_H}} = \{r \in R \mid \beta_h(r1_{h^{-1}}) = r1_h,\,\ \text{for each}\,\, h \in \mathcal{S}_H\}.$

\vu

{\bf Proof:} \em Since $R$ is a $\beta$-Galois extension of $R^{\beta}$, $R$ is a $\beta_H$-Galois extension of $R^{\beta_H}$, by Proposition \ref{propoprimeira} (note that $H$ is wide, hence $R = \bigoplus_{e \in G_0}E_e = \bigoplus_{e \in H_0}E_e = R_H$). Then, $V_R(R^{\beta_H}) = \bigoplus_{h \in H}J_h$, by Lemma \ref{comutador}. But $J_h = \{0\}$ for each $h \in H\setminus\mathcal{S}_H$, thus $\gamma(H) = \bigoplus_{h \in \mathcal{S}_H}J_h  = V_R(R^{\beta_H})$.

Noting that $C(R) \subseteq V_R(R^{\beta_H})$, we have that $\gamma(H)$ is a subalgebra of $R$ over $C(R)$. Moreover, since $\mathcal{S}_H \subseteq H$, $R^{\beta_H} \subseteq R^{\beta_{\mathcal{S}_H}}$. So $V_R(R^{\beta_{\mathcal{S}_H}}) \subseteq V_R(R^{\beta_H}) = \gamma(H)$.

\vu

Conversely, we claim that $\gamma(H) \subseteq V_R(R^{\beta_{\mathcal{S}_H}})$. Indeed, for each $r \in J_h$, with $h \in \mathcal{S}_H$, $r\beta_h(x1_{h^{-1}}) = xr$ for all $x \in R$, in particular for all $x \in R^{\beta_{\mathcal{S}_H}}$. It means that $rx = xr$ for all $x \in R^{\beta_{\mathcal{S}_H}}$, concluding that $J_h \subseteq V_R(R^{\beta_{\mathcal{S}_H}})$ for each $h \in \mathcal{S}_H$. Therefore,
$\gamma(H)=\bigoplus_{h \in \mathcal{S}_H}J_h \subseteq V_R(R^{\beta_{\mathcal{S}_H}})$.  Thus, $\gamma(H) = V_R(R^{\mathcal{S}_H})$. \hfill{$\blacksquare$}
\end{lema}

\begin{lema} By keeping the notations of Lemma \ref{igualdadegamma}, if the map $\gamma$ is injective, then there is no proper subgroupoids between $\mathcal{S}_H$ and $H$ for any wide subgropoid $H$ of $\G$. In particular, either $\mathcal{S}_H$ is a subgroupoid of $\G$ or $H = \langle\mathcal{S}_H\rangle$, the subgroupoid generated by the elements in $\mathcal{S}_H$.

\vu

{\bf Proof:} \em Suppose that there exists a proper subgroupoid $H'$ between $\mathcal{S}_H$ and $H$. Then, $R^{\beta_H} \subseteq R^{\beta_{H'}} \subseteq R^{\beta_{\mathcal{S}_H}}$ and so $V_R(R^{\beta_{\mathcal{S}_H}})\subseteq V_R(R^{\beta_{H'}}) \subseteq V_R(R^{\beta_H}) $.

By Lemma \ref{comutador}, $\gamma(H) = V_R(R^{\beta_H})$ and by Lemma \ref{igualdadegamma}, $\gamma(H) = V_R(R^{\beta_{\mathcal{S}_H}})$, thus $\gamma(H) = V_R(R^{\beta_{H'}})$. It means that $V_R(R^{\beta_{H'}}) = V_R(R^{\beta_H})$, concluding that $\gamma(H') = \gamma(H)$. But $H' \neq H$ and $\gamma$ is injective by hypothesis, which leads to a contradiction. The last assertion is obvious.\hfill{$\blacksquare$}
\end{lema}

We recall from the literature that $R$ satisfies the double centralizer property for a subring $A$ of $R$ if $V_R(V_R(A)) = A$ (\cite{demeyer}). The next theorem shows that every $C(R)$-Azumaya algebra $R$ satisfies the double centralizer property for any separable algebra $S$ over $C(R)$.

\begin{teo}\emph{\cite[Theorem 4.3]{demeyer}}\label{deme} Consider $S$ a commutative ring, $R$ a central separable $S$-algebra and suppose that $A$ is a subalgebra of $R$ which contains $S$. Then $V_R(A)$ is a separable subalgebra of $R$ and $V_R(V_R(A)) = A$. If $A$ is central, then $V_R(A)$ is also central, and the map of $S$-algebras $A \otimes V_R(A) \longrightarrow R$ given by $a \otimes r \mapsto ar$ is an isomorphism.
\end{teo}

\vu

The two theorems in the sequel are the last main theorems of this section.

\vu

\begin{teo}\label{problema} By keeping the notations of Lemma \ref{igualdadegamma}, if $R$ satisfies the double centralizer property for $R^{\beta_H}$ for each subgroupoid $H$ of $\G$, then $\theta$ is injective if and only if $\gamma$ is injective.

\vu

{\bf Proof:} \em Since $\gamma(H) = \bigoplus_{h \in H}J_h = V_R(R^{\beta_H}) = V_R(\theta(H)) = (V_R \circ \theta)(H)$ for each subgroupoid $H$ of $\G$, we have that $\gamma = V_R \circ \theta$, where in this case we are considering $V_R$ as the map which leads each subring of $R$ into its commutator.

\vu

($\Leftarrow$) Assume $\theta$ injective and $R$ satisfying the double centralizer property for $R^{\beta_H}$. Suppose $\gamma(H) = \gamma(L)$, for $H, L$ subgroupoids of $\G$. Then $(V_R \circ \theta)(H) = (V_R \circ \theta)(L)$, that is, $V_R(R^{\beta_H}) = V_R(R^{\beta_L})$. Thus, $V_R(V_R(R^{\beta_H})) = V_R(V_R(R^{\beta_L}))$. Applying the double centralizer property, we have that $R^{\beta_H} = R^{\beta_L}$, that is, $\theta(H) = \theta(L)$. Since $\theta$ is injective, we obtain $H = L$.

\vu

($\Rightarrow $) It is obvious since if $\gamma$ is injective, clearly $\theta$ so also is. \hfill{$\blacksquare$}

\end{teo}

We end this section by characterizing when two subgroupoids $H$ and $L$ of $\G$ are such that $\gamma(H) = \gamma(L)$.

\begin{teo}\label{gammatheta} Assume that, for each subgroupoid $F$ of $\G$, $R$ satisfies the double centralizer property for $R^{\beta_F}$. Then for any two subgroupoids $H$ and $L$ of $\G$, the following conditions are equivalents:
\begin{itemize}
\item[(i)] $\gamma(H) = \gamma(L)$;

\item[(ii)] $\theta(H) = \theta(L)$;

\item[(iii)] $J_g = \{0\}$ for each $g \in H \vee L\setminus H \cap L$.
\end{itemize}

\vu

{\bf Proof:} \em (i)$\Rightarrow$(ii) Let $\gamma(H) = \gamma(L)$ for subgroupoids $H, L$ of $\G$. Then $V_R(R^{\beta_H}) = \bigoplus_{h \in H}J_h = \gamma(H) = \gamma(L) =$ $\bigoplus_{l \in L}J_l = V_R(R^{\beta_L})$. Thus, $V_R(V_R(R^{\beta_H})) = V_R(V_R(R^{\beta_L}))$. Since $R$ satisfies the double centralizer property for $R^{\beta_H}$ and for $R^{\beta_L}$, we conclude that $R^{\beta_H} = R^{\beta_L}$, that is, $\theta(H) = \theta(L)$.

\vu

(ii)$\Rightarrow$(i) As it was seen in the proof of  Theorem \ref{problema}, $\gamma = V_R \circ \theta$. Thus, $\gamma(H) = V_R(\theta(H)) = V_R(\theta(L)) = \gamma(L)$.

\vu

(ii)$\Rightarrow$(iii) By assumption, $\theta(H) = \theta(L)$, that is, $R^{\beta_H} = R^{\beta_L}$. So, as it was seen in the proof of the Lemma \ref{theta}, $\bigoplus_{g \in H\vee L} J_g =$ $\bigoplus_{h \in H} J_h =$ $\bigoplus_{l \in L} J_l.$ Noting that $H \subseteq H\vee L$, we have that $\bigoplus_{g \in H\vee L} J_g =$ $\bigoplus_{h \in H} J_h$ implies that $J_g = 0$ for each $g \in H\vee L\setminus H$. Similarly,
$\bigoplus_{g \in H\vee L} J_g =$ $\bigoplus_{l \in L} J_l$ implies that $J_g = 0$ for each $g \in H\vee L\setminus L$. It means that $J_g = 0$ for each $g \in H\vee L\setminus H \cap L)$.

\vu

(iii)$\Rightarrow$(i) Since $J_g = \{0\}$ for each $g \in H\vee L\setminus H \cap L$, $\bigoplus_{g \in H \cap L} J_g = \bigoplus_{g \in H\vee L}=
\bigoplus_{h \in H} J_h = \bigoplus_{l \in L} J_l.$ Thus $\gamma(H) = \gamma(L)$.\hfill{$\blacksquare$}
\end{teo}


\section{Central $\beta$-Galois Algebras with bijective Galois map}

Throughout this section, we maintain the same assumptions on $R$, $\G$ and $\beta$, assumed in the previous one. In Section 3 we give conditions for the Galois map, corresponding to the groupoid action $\beta$, to be injective. In this section we restrict such a study to the particular case where the Galois extensions are central Galois algebras. Actually, we will prove that  the Galois map $\theta: H\mapsto R^{\beta_H}$, from the set of all the wide subgroupoids of $\G$ to the set of the subalgebras of $R$ which are separable over $R^\beta$, is indeed bijective. This result is also known as the fundamental theorem of the Galois theory for central Galois algebras.

\vu

For the sequel consider the following notations, for any subring $S$ of $R$: $$H_S=\{g\in\G\ |\ \beta_g(s1_{g^{-1}})=s1_g,\,\text{for all}\, s\in S\} \quad\text{and} \quad S'=V_R(S).$$

\vu

\begin{teo}\label{idafund} Assume that $R$ is a central $\beta$-Galois algebra over $R^{\beta}$.  If $R$ satisfies the fundamental theorem, then for any subalgebra $S$ of $R$, separable over $R^\beta$, $V_R(S) = \bigoplus_{g \in H_S}J_g$ and $S = \bigoplus_{g \in H_{S'}} J_g$.

\vd

{\bf Proof:} \em Since $S$ is a subalgebra of $R$ separable over $R^\beta$ and $R$ satisfies the fundamental theorem, then $S = R^{\beta_{H_S}}$.  Thus, $R$ is a $\beta_{H_S}$-Galois extension of $S$, by Proposition \ref{propoprimeira}, and $V_R(S) = V_R(R^{\beta_{H_S}}) = \bigoplus_{g \in H_S}J_g$, by  Lemma \ref{comutador}.

\vu

Furthermore, $R$ is an Azumaya algebra, by Remark \ref{hseparavel}. Thus, applying Theorem \ref{deme}, we have that $S' = V_R(S)$ is a subalgebra of $R$ separable over $R^{\beta}$ and $V_R(S') = V_R(V_R(S)) = S$.

Moreover, since $R$ satisfies the fundamental theorem, we also have that $S' = R^{\beta_{H_{S'}}}$. Therefore, $R$ is a $\beta_{H_{S'}}$-Galois extension of $S'$, again by Proposition \ref{propoprimeira}, and, consequently, $S = V_R(S') = V_R(R^{\beta_{H_{S'}}}) = \bigoplus_{g \in H_{S'}} J_g$, by Lemma \ref{comutador}.\hfill{$\blacksquare$}

\begin{lema}\label{sep} Let $H$ be a wide subgroupoid of $\G$. If $R$ is separable algebra over $R^{\beta}$, then $R$ is separable over $R^{\beta_H}$.

{\bf Proof:} \em Consider the map
$$
\begin{array}{cccc}
\varphi \ : & \! R \times R & \! \longrightarrow
& \! R \otimes_{R^{\beta_H}} R \\
& \! (r, s) & \! \longmapsto
& \! r \otimes_{R^{\beta_H}} s
\end{array}
$$
It is straightforward to check that $\varphi$ is additive and $R^{\beta}$-balanced, hence it induces a map
$$
\begin{array}{cccc}
\overline{\varphi} \ : & \! R \otimes_{R^{\beta}} R & \! \longrightarrow
& \! R \otimes_{R^{\beta_H}} R, \\
& \! r \otimes_{R^{\beta}} s & \! \longmapsto
& \! r \otimes_{R^{\beta_H}} s
\end{array}
$$ which also is a homorphism of $R$-bimodules.

\vu

By hypothesis, $R$ is separable over $R^{\beta}$, thus there exist $e = \sum_i a_i \otimes_{R^{\beta}} b_i$ $\in R \otimes_{R^{\beta}} R$ such that $er = re$ for all $r \in R$ and $\sum a_ib_i = 1_{R}$. Consider $e' = \overline{\varphi}(e) = \sum_i a_i \otimes_{R^{\beta_H}} b_i \in R \otimes_{R^{\beta_H}} R$. Thereby, $$e'r = \overline{\varphi}(e)r = \overline{\varphi}(er) = \overline{\varphi}(re) = r\overline{\varphi}(e) = re',$$ for all $r \in R$, and $\sum a_ib_i = 1_{R}$, showing that $R$ is separable over $R^{\beta_H}$.\hfill{$\blacksquare$}

\end{lema}

The next theorem is the converse of Theorem \ref{idafund}.

\begin{teo}\label{voltafund} Suppose that $R$ is a central $\beta$-Galois algebra separable over $R^{\beta}$. If for any subalgebra $S$ of $R$ separable over $R^{\beta_H}$, $V_R(S) = \bigoplus_{g \in H_S}J_g$, then $R$ satisfies the fundamental theorem.

{\bf Proof:} \em Since $R$ is a central $\beta$-Galois algebra, the map $\theta: H \mapsto R^{\beta_H}$ is injective, by Corollary \ref{hinjetiva}. Furthermore, since $R$ is separable over $R^{\beta}$, it follows that $R$ is separable over $R^{\beta_H}$, by Lemma \ref{sep}.

Now let $S$ be a subalgebra of $R$ separable over $R^{\beta}$. Then, by Theorem \ref{deme}, $S' = V_R(S)$ is also a subalgebra of $R$ separable over $R^{\beta}$ and $V_R(S') = V_R(V_R(S)) = S$.

On the other hand, $R$ is a $\beta_{H_S}$-Galois extension of $R^{\beta_{H_S}}$, by Proposition \ref{propoprimeira}, thus $V_R(R^{\beta_{H_S}}) = \bigoplus_{g \in H_S}J_g$, by Lemma \ref{comutador}. So, $V_R(S) = \bigoplus_{g \in H_S}J_g = V_R(R^{\beta_{H_S}})$. Moreover, since $H_S$ is a wide subgroupoid of $\G$, $R$ is separable over $R^{\beta_{H_S}}$, by the Lemma \ref{sep}. Therefore, by the Theorem \ref{deme},$$S = V_R(V_R(S)) = V_R(V_R(R^{\beta_{H_S}})) = R^{\beta_{H_S}}.$$ It means that the map $\theta: H \mapsto R^{\beta_H}$ is surjective, which completes the proof.\hfill{$\blacksquare$}
\end{teo}

To end this section , we summarize the above results in the following theorem, which gives a characterization of central $\beta$-Galois algebras that satisfy the fundamental theorem.

\begin{teo} \label{teolast} Assume that $R$ is a central $\beta$-Galois algebra over $R^{\beta}$. Then $R$ satisfies the fundamental theorem if and only if $V_R(S) = \bigoplus_{g \in H_S}J_g$, for any subalgebra $S$ of $R$ separable over $R^{\beta}$.

\end{teo}

\end{teo}

\bibliographystyle{amsalpha}
{
}

\end{document}